\magnification=1200

{\bf Positivity conditions for bihomogeneous polynomials}

David W. Catlin and John P. D'Angelo
\bigskip

{\bf Introduction}

In this paper we continue our study of
a complex variables version of Hilbert's seventeenth problem
by generalizing some of the results from [CD]. 
Given a bihomogeneous polynomial $f$
of several complex variables that is positive away from the
origin, we proved that there is an integer $d$ so that $||z||^{2d} f(z,{\overline z})$
is the squared norm of a holomorphic mapping. Thus, although $f$ may not itself
be a squared norm, it must be the quotient of squared norms of holomorphic homogeneous
polynomial mappings. The proof required some
operator theory on the unit ball. In the present paper we
prove that we can replace the squared Euclidean norm by squared norms 
arising from
an orthonormal basis for the space of homogeneous polynomials on any 
bounded circled pseudoconvex domain of finite type. To do so we prove
a compactness result for an integral operator on such domains related to the
Bergman kernel function. Recall that the Bergman kernel function $B$ for
a domain $\Omega$ is the integral
kernel for the operator $P$ that projects $L^2(\Omega)$ to the closed subspace
$A^2(\Omega)$ of holomorphic functions in $L^2(\Omega)$.

We prove the following results.

{Proposition 1}. Suppose that $\Omega$ is a bounded pseudoconvex domain in
${\bf C^n}$ for which the
${\overline \partial}$-Neumann operator $N$ is compact. Let $M$ be a pseudodifferential
operator of order 0. Then the commutator $[P,M]$ is compact on $L^2(\Omega)$.

{Theorem 1}. Suppose that $\Omega$ is a smoothly
bounded pseudoconvex domain of finite type
in ${\bf C^n}$,  with Bergman kernel $B(z,{\overline \zeta})$.
Let $g$ be a smooth function on ${\overline \Omega} \times {\overline
\Omega}$ that vanishes on the boundary diagonal. 
Then the operator on $L^2(\Omega)$ with integral kernel 
$gB $ is compact.

{Theorem 2}. Suppose that $\Omega$ is a smoothly bounded pseudoconvex circled domain 
in ${\bf C^n}$ of finite type. For each integer $d$, let
$\Phi^d = (\Phi^d_1,...,\Phi^d_N) $ denote an orthonormal basis for the homogeneous
polynomials of degree $d$ on $\Omega$. Let $f$ be a bihomogeneous polynomial that
is positive away from the origin. Then there is an integer $d_0$ (depending on $f$)
such that, for each $d \ge d_0$, there is a homogeneous polynomial
mapping $h$ such that
$$ ||\Phi^d(z)||^2 f(z,{\overline z}) = ||h(z)||^2 $$

In section III we interpret 
Theorem 2 and the theorem from [CD] in terms of Hermitian line bundles over
complex projective space. In a future paper we will
prove a related differential geometric theorem involving metrics on Hermitian
line bundles over compact complex manifolds.

The second author acknowledges useful discussions with Yum-Tong Siu, Steve Bradlow,
and Alex Tumanov. He also acknowledges support from MSRI. The first author
acknowledges support from NSF. The authors acknowledge helpful comments by the referee.

\bigskip

{\bf I. Needed facts about compact operators}.

In this section we write
$(u,v)$ for the inner product on a Hilbert space, and $||u||^2$ for
the squared norm. An operator 
$A:H_1 \to H_2$ between Hilbert spaces is compact
if whenever $f_j$ is a bounded sequence in $H_1$, then $A(f_j)$ has a convergent
subsequence in $H_2$. We convert this into an estimate:

{\bf Lemma 1}. Suppose that $A:H_1 \to H_2$ is an operator between Hilbert spaces.
Then $A$ is compact if and only if, 
for all positive $\epsilon$,  there is a compact operator $B=B_\epsilon$ and
a positive constant $C_\epsilon$ such that
the following estimate holds:

$$ ||Af||^2 \le \epsilon ||f||^2 + C_\epsilon ||B_\epsilon f||^2 \eqno (1) $$

Proof. If $A$ is compact, we may take $B=A$ and $C_\epsilon = 1$. Conversely,
we suppose that (1) holds. Let $f_j$ be a bounded sequence in the domain; 
we will construct a subsequence whose image under $A$ is Cauchy. 
Consider a positive integer $n$.
From (1) we have

$$ ||A(f_j) - A(f_k)||^2 \le \epsilon ||f_j -f_k||^2 +
 C_\epsilon ||B_\epsilon(f_j- f_k)||^2 . \eqno (2) $$

Since $||f_j||$ are bounded, the first term 
can be made smaller than ${1 \over 2 n}$ by choosing
$\epsilon$ small enough. The second can then be made smaller than ${1 \over 2n}$
 by extracting a subsequence (still labeled the same)
for which $B_\epsilon f_j$ converges and then choosing
$j$ and $k$ sufficiently large. Thus for each $n$ there is a subsequence 
such that $||A(f_{j(n)}) - A(f_{k(n)})||^2 < {1 \over n}$.
Extracting the diagonal subsequence gives a subsequence
whose image under $A$ is Cauchy. Therefore $A$ is compact.
$\spadesuit$

In order to prove (1) 
in specific cases we will use the standard remark that, given $\epsilon > 0$,
there is a  positive constant $C_\epsilon$ so that
$$ |(u,v)| \le \epsilon ||u||^2 + C_\epsilon ||v||^2 .\eqno (3) $$
We sometimes write $sc$ for $\epsilon$ and $lc$ for $C_\epsilon$.

Let $\Omega$ be a smoothly bounded pseudoconvex domain in ${\bf C^n}$.
We assume that the reader is familar with the ${\overline \partial}$-Neumann
problem. See [C2, D,FK,K] for example.
We have the Hodge decomposition on $(0,q)$ forms
given by 
$ I =  ({\overline \partial}^* {\overline \partial} + 
{\overline \partial} {\overline \partial}^*) N + H $
where $H$ is the harmonic projector. For pseudoconvex
domains in ${\bf C^n}$, the operator $H$ on $(0,q)$ forms vanishes except
when $q=0$, in which case it equals the Bergman projection $P$.
Thus on $(0,q)$ forms for $q \ge 1$, the ${\overline
\partial}$-Neumann operator $N$ satisfies
$({\overline \partial}^* {\overline \partial} 
+ {\overline \partial} {\overline \partial}^* ) N = I$.

The Bergman kernel function $B(z,{\overline w})$
for a bounded domain $\Omega$ is the integral
kernel of the operator $P$
that projects $L^2(\Omega)$ onto the closed subspace $A^2(\Omega)$
of  holomorphic functions in
$L^2(\Omega)$. 
Suppose that the collection $\{\phi_\alpha\}$ forms a complete orthonormal set for
$A^2(\Omega)$. Then the Bergman kernel satisfies 
$$ B(z,{\overline \zeta}) = \sum \phi_\alpha(z) {\overline {\phi_\alpha(\zeta)}} .
\eqno (4) $$

Kohn's formula relates
the $N$ operator to the Bergman projection $P$; 
it states that
$P = I - {\overline \partial}^* N {\overline \partial}$.
Here the $N$ operator is defined on $(0,1)$ forms.
The image of $N$ is contained in the domain of ${\overline \partial}^*$.
For forms of all degrees this remains true; then 
${\overline \partial}^* N$ maps $(0,q)$ forms in $L^2(\Omega)$
to $(0,q-1)$ forms in $L^2$, and its adjoint 
 $N {\overline \partial}$ maps $(0,q-1)$ forms in $L^2$ to $(0,q)$ forms in $L^2$.

We recall (See [D]) that $\Omega$ is of {\it finite type} if there is a bound on the
order of contact of ambient complex analytic varieties with ${\rm b}\Omega$ at all
points. Domains of finite type satisfy subelliptic estimates [C1], and hence the
${\overline \partial}$-Neumann 
operator is compact. Another consequence [Ke] is that the Bergman kernel
function is smooth on ${\overline \Omega} \times {\overline \Omega}$ away from the
boundary diagonal. See [HI] for recent results concerning compactness of the
${\overline \partial}$-Neumann operator.

{\bf Lemma 2}. Suppose that the ${\overline \partial}$-Neumann operator $N$
is compact on the spaces of forms of type $(0,q)$ and $(0,q+1)$ in $L^2(\Omega)$.
Then ${\overline \partial}^* N$ and $N {\overline \partial}$ are compact.

Proof. Since $N {\overline \partial}$ on $(0,q)$ forms is the adjoint of 
${\overline \partial}^* N$ on $(0,q+1)$ forms, it suffices to prove that 
${\overline \partial}^* N$
is compact. Using the Hodge decomposition formula 
$ I =  
({\overline \partial}^* {\overline \partial} + 
{\overline \partial} {\overline \partial}^*) N + H $
we see that

$$ || {\overline \partial}^* Nf||^2 = 
({\overline \partial} {\overline \partial}^* Nf, Nf) =
(f,Nf) - ({\overline \partial}^* {\overline \partial} N f, Nf) - (Hf,Nf) . \eqno (5) $$
Since $(Hf,Nf) = 0$, and 
$- ({\overline \partial}^* {\overline \partial} N f, Nf) = 
- ||   {\overline \partial}Nf||^2 $,  we see that

$$ || {\overline \partial}^* Nf||^2 \le (f,Nf) \le \epsilon ||f||^2 + C_\epsilon
||Nf||^2 $$

Lemma 1 now implies that ${\overline \partial}^* N $ is compact.
$\spadesuit$.

The next proposition is a key step in our proof of Theorem 1. It also
could be used 
instead of Theorem 1 in the proof of Theorem 2.

{\bf Proposition 1}. Suppose that $\Omega$ is a bounded pseudoconvex 
domain in ${\bf C^n}$ for which the
${\overline \partial}$-Neumann operator $N$ 
is compact on $(0,1)$ forms in $L^2(\Omega)$. Let $M$ be a pseudodifferential
operator of order 0. Then the commutator $A = [P,M]$ is compact on $L^2(\Omega)$.

Proof. By the formula relating $N$ and $P$ we have 

$$ A = [P,M] = [I - {\overline \partial}^* N {\overline \partial}, M] = 
 -[{\overline \partial}^* N {\overline \partial}, M] 
= - [{\overline \partial}^*,M] N {\overline \partial} + 
{\overline \partial}^* [M, N {\overline \partial} ] \eqno (6) $$

Note first that the commutator of $M$ and either $\overline \partial$ or
${\overline \partial}^*$ is an operator of order zero, and hence bounded.
The first term in the last equality in (6) is the composition of the bounded operator
$ [{\overline \partial}^*,M]$ with the compact operator 
$N {\overline \partial}$, and hence is itself compact. The second term is more
difficult. We use both lemmas.
As usual we write 

$$ Q(u,v) = ({\overline \partial}u,{\overline \partial}v) + 
({\overline \partial}^*u,{\overline \partial}^*v) $$

The main property of $Q$ is that  $Q(Nu,v) = (u,v)$.
We write $a = MN{\overline \partial}f - N{\overline \partial}Mf $ to simplify notation.
Our goal is to show that the map taking $f$ to ${\overline \partial}^* a$ is compact.
To do so, we compute
$Q(a,a)$.

When we commute first order operators past $M$ we obtain operators of order zero,
all denoted by $P_0$. We obtain

$$ Q(MN{\overline \partial}f-N{\overline \partial}Mf,
MN{\overline \partial}f-N{\overline \partial}Mf ) = 
||{\overline \partial}(MN {\overline \partial}f-N{\overline \partial}M)f||^2 
+ ||{\overline \partial}^* (MN {\overline \partial}f-N{\overline \partial}M)f||^2 
$$
$$ = ({\overline \partial}^* MN{\overline \partial}f,{\overline \partial}^* a)
- ({\overline \partial}^* N{\overline \partial}Mf,{\overline \partial}^* a)
+  ({\overline \partial} MN{\overline \partial}f,{\overline \partial} a)
- ({\overline \partial} N{\overline \partial}Mf,{\overline \partial} a) \eqno (7)
$$
Using $Q(Nu,v) = (u,v)$, the second and fourth terms in (7) simplify to 
$- ( {\overline \partial}Mf,a)$. Using this, commuting $M$ past the
differentiations, and moving it to the other side, we obtain
$$
= (P_0 N {\overline \partial}f,  {\overline \partial}^* a)
+ ({\overline \partial}^* N{\overline \partial}f, M^* {\overline \partial}^*a)
- ( {\overline \partial}Mf,a) 
+ (P_0 N {\overline \partial}f, {\overline \partial} a)
+ ({\overline \partial} N {\overline \partial}f,  M^* {\overline \partial} a) $$

$$ = (P_0 N {\overline \partial}f,  {\overline \partial}^* a)
+ ({\overline \partial}^* N{\overline \partial}f, {\overline \partial}^* M^* a)
+ ({\overline \partial}^* N{\overline \partial}f, P_0 a)
- ({\overline \partial}Mf,a) $$
$$ 
+   (P_0 N {\overline \partial}f,  {\overline \partial} a)
+  ({\overline \partial} N{\overline \partial}f, {\overline \partial}M^* a)
+ ({\overline \partial} N{\overline \partial}f, P_0 a) \eqno (8) $$
Using $Q(Nu,v) = (u,v)$, the second term plus the sixth term in (8) becomes
$({\overline \partial}f,M^*a)$.  Subtracting the fourth term from this
gives $([M,{\overline \partial}]f,a)$. Hence we obtain

$$ Q(a,a) = ( P_0 N {\overline \partial}f,  {\overline \partial}^* a) +
(P_0 f,a) +  ({\overline \partial}^* N{\overline \partial}f, P_0 a) 
+  ( P_0 N {\overline \partial}f,  {\overline \partial} a)
+ ({\overline \partial} N{\overline \partial}f, P_0 a)
$$
$$ = T_1 + T_2 + T_3 + T_4 + T_5 \eqno (9) $$

Our next goal is to show that 
$$ Q(a,a) \le sc ||f||^2 + lc ||Bf||^2 $$
where $B$ is compact.
We estimate each $T_j$:
$$ T_1 \le sc ||{\overline \partial}^* a||^2 + lc ||N{\overline \partial}f||^2 .$$

$$ T_2 \le sc ||f||^2 + lc ||a||^2 .$$
We have $$ T_3 \le sc ||{\overline \partial}^* N {\overline \partial} f||^2 + 
lc ||a||^2 = sc ||(I-P)f||^2 + lc ||a||^2 \le sc ||f||^2 + lc ||a||^2 .$$

$$ T_4 \le sc ||{\overline \partial} a||^2 + lc || N{\overline \partial }f||^2 .$$

We claim that ${\overline \partial} N {\overline \partial} f =
0$, so $T_5$ vanishes. This follows because the ${\overline \partial}$ and $N$
operators, defined on forms of different degrees, commute. Another proof writes 
$ f = Hf + 
({\overline \partial} {\overline \partial}^* + 
{\overline \partial}^* {\overline \partial}) Nf$,
so 
${\overline \partial}f = 
{\overline \partial} {\overline \partial}^* {\overline \partial} Nf$, and thus 
$$ {\overline \partial} N {\overline \partial} f = 
{\overline \partial} N 
{\overline \partial} {\overline \partial}^* {\overline \partial} Nf=
{\overline \partial} N ({\overline \partial} {\overline \partial}^* +
{\overline \partial}^* {\overline \partial} ) {\overline \partial} N f =
{\overline \partial} {\overline \partial} N f = 0 .$$

Now we can subtract terms of the form  $sc ||{\overline \partial}a||^2$ and
$sc ||{\overline \partial}^* a||^2$ and obtain an estimate for $Q(a,a)$ involving the
other terms. By Lemma 2, $N{\overline \partial}$ is compact.
Composing it with $M$ on either side still gives a compact operator,
because $M$ is bounded. Hence the map sending
$f$ to $a = [M,N{\overline \partial}]f$ is compact. Thus we can absorb all the terms
on the right in terms of $C ||Bf||^2$, where $B$ is compact, except for the $sc
||f||^2$ terms. Since $||{\overline
\partial}^*a||^2 \le Q(a,a)$, we obtain

$$ ||{\overline \partial }^* a||^2 \le Q(a,a) \le 
sc ||f||^2 + C ||Bf||^2 . \eqno (10) $$
Statement (10) and Lemma 1 imply the desired conclusion. $\spadesuit$

We will use Theorem 1 in order to prove Theorem 2; 
a weaker version (for polynomials vanishing on the diagonal) suffices for
our application. If we assume that $g$ is a polynomial vanishing on the diagonal,
then we can prove Theorem 1 assuming only that $N$ is compact by using Proposition 1.

{\bf Theorem 1}. Suppose that $\Omega$ is a smoothly bounded pseudoconvex
domain of finite type with Bergman kernel function $B$. Suppose that
$g$ is smooth on ${\overline \Omega} \times {\overline \Omega}$ and that $g$
vanishes on the boundary diagonal. 
Let $T:L^2(\Omega) \to L^2(\Omega)$ be the operator whose kernel is given by
$gB$.
Then $T$ is a compact operator.

Proof.  We will write $g(z,\zeta)$ for the value of $g$;
this notation does not mean that $g$ is holomorphic in either variable!
First we write $z= \zeta + (z-\zeta)$.
We can therefore write by Taylor's formula 

$$ g(z, \zeta) = g(\zeta +(z- \zeta),\zeta) = g_0(\zeta) + 
\sum_{|\alpha|,|\beta| \le N}
f_{\alpha \beta}(\zeta) (z- \zeta)^\alpha ({\overline
z}-{\overline \zeta})^\beta +R_N \eqno (11)
$$

Here $R_N$ is the remainder term.
Evaluating (11) at $z = \zeta$ gives $g(\zeta,\zeta) = g_0(\zeta)$. Since
$g$ vanishes on the boundary diagonal we conclude that 
$g_0$ vanishes on the boundary.

Each of these terms gets multiplied by
the Bergman kernel. We first
show that $g_0(\zeta) B(z,{\overline \zeta})$ defines a compact
operator ${\cal A}$.
To establish the compactness, we use Lemma 1. Suppose that $r$ is a defining function
for $\Omega$.
Given $\epsilon >0$, we can find $\delta > 0$ so that
$-\delta \le r \le 0$ implies $|g_0| \le \epsilon$.
Choose a smooth positive function $\chi$, bounded by unity, that equals unity when $r
<-\delta$, and is supported in $\Omega$. Then we write

$$ g_0 B = \chi g_0 B + (1 - \chi) g_0 B = A_1 + A_2 $$

The operator ${\cal A}_1$
defined by $A_1$ is compact, since its kernel is smooth on all of
${\overline \Omega} \times {\overline \Omega} $. 
The operator ${\cal A}_2$ defined by $A_2$ is multiplication by
$g_0 (1 - \chi)$ followed by the Bergman projection $P$.
Note that $|g_0| \le \epsilon$ where $\chi$ is not equal to one.
Hence we can write

$$ ||{\cal A} f || = || {\cal A}_1 f + {\cal A}_2 f|| 
\le ||{\cal A}_1 f|| + || P(1-\chi)g_0 f|| 
\le  ||{\cal A}_1 f|| + ||(1-\chi)g_0 f|| 
\le ||{\cal A}_1 f|| + \epsilon ||f||$$ 

Thus ${\cal A}$ defines a compact operator by Lemma 1.

We now consider the sum in (11). Each of the terms 
$ f_{\alpha \beta}(\zeta) (z- \zeta)^\alpha ({\overline
z}-{\overline \zeta})^\beta $
is a finite sum of terms in the form
$ A_j(z) B_j(\zeta) (z_j - \zeta_j)$ or 
$ A_j(z) B_j(\zeta) ({\overline z_j }- {\overline \zeta_j})$
where $A_j$ and $B_j$ are smooth. The composition in either order
of a compact operator with a bounded operator
is a compact operator, and a finite sum of compact operators is a compact
operator. Therefore, to establish the compactness  of $A$, it suffices to prove
the following:
Suppose that $A:L^2 \to L^2$ is the operator whose kernel
is given by $ ({\overline \zeta_j} - {\overline z_j}) B(z, {\overline \zeta})$, or
by  $(z_j - \zeta_j) B(z,{\overline \zeta})$.
Then $A$ is compact.

Thus the crucial point is to show that such an $A$ is compact. 
Note however that such an
$A$ is the commutator $[P,M]$, where $M$ denotes multiplication
by ${\overline z_j}$ or its conjugate. Thus the compactness of $A$ follows
from Proposition 1.

Finally we must handle the remainder term from (11). It suffices to show that
the kernel $R_N(z,\zeta) B(z,{\overline \zeta})$ is continuous on
${\overline \Omega} \times {\overline \Omega}$. To see this it suffices to
show that there is a positive exponent $\alpha$ such that 

$$ |B(z,{\overline \zeta})| \le C |z-\zeta|^{-\alpha}. \eqno (12) $$
If (12) holds, then we can choose the index $N$ in the Taylor expansion (11)
so large  that
 $R_N(z,\zeta) B(z,{\overline \zeta})$ is continuous on
${\overline \Omega} \times {\overline \Omega}$, and therefore defines a compact
operator. We state and prove (12)  in 
Lemma 3; its proof completes the proof of Theorem 1.
$\spadesuit$.

{\bf Lemma 3}. Suppose that $\Omega$ is a smoothly bounded pseudoconvex domain for which
a subelliptic estimate of order $\epsilon$ holds on $(0,1)$ forms. Then there is a
constant $C$ so that 

$$ |B(z,{\overline \zeta})| \le C |z-\zeta|^{-{(2n+6) \over \epsilon}} \eqno (13) $$

Proof. The proof follows the methods of [Ke] and [C2]. 
First we assume that $z \in {\overline
\Omega}$, that $\zeta \in \Omega$, and that the distance between them is $t$.
Choose $\delta $ equal to the minimum of ${t \over 4}$ and the distance between
$\zeta$ and the boundary. 
Let $\phi_\zeta$ be an approximation to the identity at $\zeta$, spherically
symmetric and supported in a ball about radius $\delta$ about $\zeta$. We may
write

$$ B(z,\zeta) = (I - {\overline \partial}^* N {\overline \partial} )(\phi_\zeta)(z) $$

Choose smooth cut-off functions $\chi_1$ and $\chi_2$ supported in a ball of radius
${t \over 4}$ and so that $\chi_2 = 1$ near
$supp(\chi_1)$. Using the methods of Proposition 2 from
[C2], there is an estimate involving Sobolev
norms:

$$ |\chi_1 B(.,\zeta)| \le ||\chi_1 {\overline \partial}^* N {\overline
\partial}\phi_\zeta||_{n+1} \le C t^{- {n+3 \over \epsilon}} || \chi_2 N {\overline
\partial \phi_\zeta}|| \eqno (14) $$

We estimate $ || \chi_2 N {\overline
\partial \phi_\zeta}|| $ by pairing $ \chi_2 N {\overline
\partial \phi_\zeta}$ with a smooth $(0,1)$ form $f$ with $||f||_0 \le 1$. 
We obtain using
integration by parts that 

$$ |( \chi_2 N {\overline
\partial \phi_\zeta},f)| = |(\phi_\zeta, {\overline \partial}^* N \chi_2 f)|\le
 |(\phi_\zeta, \lambda {\overline \partial}^* N \chi_2 f)| \eqno (15) $$
where we have written $\lambda$ for a cut-off function that is
unity near $supp(\phi_\zeta)$, vanishes near $supp(\chi_2)$, and whose $k$-th derivative
can be estimated by $C_k t^{-k}$. Using the generalized Schwartz inequality we can
estimate (15) by 

$$ (15) \le ||\phi_\zeta ||_{-n-1} \ 
|| \lambda {\overline \partial}^* N \chi_2 f||_{n+1} \le C
||\phi_\zeta||_{-n-1} t^{-{n+3 \over \epsilon}} ||f|| \eqno (16) $$

Taking the supremum over $f$ and estimating $||\phi_\zeta||_{-n-1}$ by a 
positive constant, we obtain

$$ |B(z,\zeta)| \le C t^{-{(2n+6) \over \epsilon}}, $$
which gives (13).
Finally, since $B$ is smooth on ${\overline \Omega} \times {\overline \Omega}$ off the
diagonal, the same estimate holds when $z, \zeta \in {\overline \Omega}$.
$\spadesuit$

\bigskip

{\bf II. Proof of the Main result}

In this section we write $||z||^2$ for the 
Euclidean norm of a point in ${\bf C^n}$.
Let $\Omega$ be a bounded domain in ${\bf C^n}$. We now write
$\langle,\rangle_\Omega$ for the $L^2$ inner product given by

$$ \langle g,h \rangle_\Omega = \int_\Omega g(z) {\overline h(z)} dV(z) $$
and write $|| g||_\Omega^2$ for $\langle g,g \rangle_\Omega $.
Recall that $\Omega$
is {\it circled} if it is invariant under scalar multiplication by
$e^{i \theta}$. We have the following simple lemma.

{\bf Lemma 4}. Suppose that $\Omega$ is a bounded circled domain. Then, for $j \ne k$,
the spaces $V_j$ and $V_k$ of
holomorphic homogeneous polynomials are orthogonal in $A^2(\Omega)$.

Proof. Suppose that $p_j$ and $p_k$ are homogeneous polynomials of the degree
indicated.
Then we have $p_j(e^{i\theta }z) = e^{i j \theta}p_j(z)$. From the definition of the
inner product as an integral, the change of variables formula for integrals, and
the invariance of $\Omega$ under multiplication by $e^{i \theta}$,  we
obtain

$$ \langle p_j, p_k \rangle_\Omega = \int_\Omega p_j(z) {\overline p_k(z)} dV(z) =
\int_\Omega  p_j(e^{i \theta z}) {\overline p_k}(e^{i \theta z}) dV(e^{i \theta z}) =
 e^{i \theta (j-k)}  \langle p_j, p_k \rangle_\Omega \eqno (17) $$
From (17) we see that the inner product must vanish unless $j=k$. $\spadesuit$

Let $N=N(d,n+1)$ denote the dimension of $V_d$.
 Suppose that $(E_{\mu \nu})$ is a
Hermitian matrix on ${\bf C^N}$. 
We consider the integral operator ${\cal E}$ defined on $V_d$ by

$$ ({\cal E} g) (z) = 
\int_\Omega \sum_{\mu,\nu} E_{\mu \nu} z^{\mu} {\overline w}^\nu g(w) dV(w) \eqno
(18) $$

{\bf Proposition 2}. Let $\Omega$ be a bounded circled domain.
The following are equivalent:

1. The Hermitian matrix $(E_{\mu \nu})$ is positive definite.

2. The operator ${\cal E}$ defined by (18) is positive definite.

3. We can write 
$$ f(z,{\overline z}) = \sum E_{\mu \nu} z^\mu {\overline z}^\nu = ||A(z)||^2 $$ 
where the components of $A$ form a basis for $V_d$. 

Proof.
Suppose that 1) holds. By linear algebra we can find basis vectors
$E_\mu$ of ${\bf C^N}$ so that $E_{\mu \nu} = \langle E_\mu, E_ \nu \rangle $.
This implies that $f(z,{\overline z}) = ||\sum E_\mu z^\mu||^2 $, so 1) implies 3).

If 3) holds we have 
$f(z,{\overline z}) = ||A(z)||^2$, where $A_1,A_2,...,A_N$ form a basis for $V_d$. Then

$$\langle {\cal E}h,h \rangle_\Omega = 
\int_\Omega \int_\Omega \sum A_j(z) {\overline {A_j(\zeta)}} h(\zeta)
{\overline {h(z)}} dV(\zeta) dV(z) = \sum |\langle A_j,h \rangle_\Omega|^2 \eqno (19) $$
and we see that ${\cal E}$ is a positive operator. Thus 3) implies 2). 

We finish the proof by showing that 2) implies 1).  
Thus we want to find $k>0$ so that

$$ \sum E_{\mu \nu} \zeta_\mu {\overline \zeta_\nu} \ge k \sum |\zeta_\mu|^2 $$
Define $ c_{\alpha \beta} = \langle z^\alpha, z^\beta \rangle_\Omega $.
Since this matrix is invertible
we may define $g_{\beta}$ by $\sum c_{\beta \nu}g_{\beta} = \zeta_\nu$.
Let $ g = \sum g_\alpha z^\alpha $. Because 
 ${\cal E}$ is positive definite, there is a $c>0$ so that 

$$ \langle {\cal E} g, g \rangle_\Omega \ge c ||g||_\Omega^2.  \eqno (20) $$
Doing
simple computations shows that 

$$ \langle {\cal E} g, g \rangle_\Omega = 
\sum E_{\mu \nu} c_{\mu \alpha} c_{\beta \nu} {\overline g_\alpha } g_{\beta} $$
$$ ||g||_\Omega^2 = \sum c_{\alpha \beta}g_\alpha {\overline g_\beta}. $$
The matrix $(c_{\alpha \beta})$ arises from inner products and hence is positive
definite. Hence there are constants
so that 
$||g||^2_\Omega \ge c' \sum |g_\alpha|^2 \ge c'' \sum |\zeta_\alpha|^2 $. 
Therefore we have 
$$ \sum E_{\mu \nu} \zeta_\mu {\overline \zeta_\nu} = 
\sum E_{\mu \nu} c_{\mu \alpha} c_{\beta \nu} {\overline g_\alpha } g_{\beta} = 
\langle {\cal E} g, g \rangle_\Omega 
\ge c||g||^2_\Omega 
\ge c c'' \sum |\zeta_j|^2 $$
This shows that 2) implies 1).
$\spadesuit$

We can now state and prove our main application of Theorem 1.

{\bf Theorem 2}. Suppose that $\Omega$ is a smoothly
bounded circled pseudoconvex domain of finite type in ${\bf C^n}$.
For each integer $d \ge 0$, let
$\Phi^d = (\Phi^d_1,...,\Phi^d_N) $ denote an orthonormal basis for the homogeneous
polynomials of degree $d$ on $\Omega$. Let $f$
be a bihomogeneous polynomial that is
positive away from the origin. Then there is an integer $d_0$ 
such that, for each $d \ge d_0$, there is a homogeneous polynomial
mapping $h$ such that
$$ ||\Phi^d(z)||^2 f(z,{\overline z}) = ||h(z)||^2 \eqno (21) $$

Proof. We will prove that, for all sufficiently large $d$,
  there is a homogenoeus
polynomial mapping $h$ such that (21) holds.
To do this we use Proposition 2; (21) holds if and only if
the operator $K$ with integral kernel

$$ f(z,{\overline \zeta}) \langle \Phi^d(z), \Phi^d(\zeta) \rangle $$
is positive on the space $V_{m+d}$.

To prove this, we first let $\xi$ be a smooth function with compact support in
$\Omega$ that is positive at the origin. We write 

$$ f(z,{\overline \zeta}) B(z, {\overline \zeta}) = 
( f(z,{\overline z}) + \xi(\zeta,{\overline \zeta}) ) B(z, {\overline \zeta}) 
 - \xi(\zeta,{\overline \zeta})  B(z, {\overline \zeta}) +
( f(z,{\overline \zeta}) - f(z,{\overline z}))  B(z, {\overline \zeta}) $$
$$ = T_1 + T_2 + T_3 \eqno (22) $$

We claim that the first term $T_1$
defines a positive operator $Q$
on all of $A^2(\Omega)$, and that
the operators defined by $T_2$ and $T_3$ are compact.
The first follows because the Bergman kernel is a self-adjoint projection.
To use this,  let $h$ be in
$A^2(\Omega)$. Then $ Qh = M_f Ph + P M_\xi h = M_f h + P M_\xi h $.
Here  $M_q$ is the operator given by multiplication by
$q$. Therefore we have 

$$ \langle Qh,h \rangle_\Omega =
\langle M_f h + P M_\xi h,h \rangle_\Omega = 
 \langle M_f h, h \rangle_\Omega +  \langle M_\xi h, Ph \rangle
= \langle M_{f+\xi} h,h \rangle_\Omega 
 \ge c ||h||^2_\Omega \eqno (23) $$

In (23) we have estimated $f+\xi$ from below by a positive constant.
This proves that $Q$ is positive.

The operator defined by $T_2$ is compact because $T_2$ is smooth on all of ${\overline
\Omega} \times {\overline \Omega}$. Since $B$ is holomorphic in its first variable,
and anti-holomorphic in its second variables, it is smooth on $\Omega \times \Omega$.
Since $\xi$ has compact
support,  $T_2$ is smooth everywhere. Given our hypotheses on $\Omega$, 
Theorem 1 implies that the  operator defined by $T_3$ is also compact.

Observe that $K$ is now known to be the sum of a compact operator $T$
and an operator $Q$ that it positive on $A^2(\Omega)$.
 Suppose that $\langle Qh, h \rangle_{\Omega} \ge c ||h||^2_\Omega $ for all $h \in
A^2(\Omega)$.
Write $|||T|||$ for the operator norm of $T$. 
Since $T$ is compact, we can find [Ru] a finite rank operator $L$ so that
$|||T-L||| < {c \over 3}$. Since $L$ is finite rank, for sufficiently large $j$ we may
assume that the restriction $L'$ of $L$ to $V_j$ has operator norm $|||L'||| < {c
\over 3}$. Therefore, on $V_j$, for sufficiently large $j$ 

$$ \langle Kh,h\rangle_\Omega = \langle Qh,h \rangle_\Omega 
+   \langle (T-L)h,h \rangle_\Omega  +  \langle Lh,h \rangle_\Omega $$
$$
\ge c||h||^2_\Omega - {c \over 3} ||h||^2_\Omega - {c \over 3} ||h||^2_\Omega = {c
\over 3} ||h||^2_\Omega \eqno (24) $$

By (24) $K$ is positive on $V_j$.
We choose a complete orthonormal set for $A^2(\Omega)$ that includes all components 
of $\Phi_d$, and replace $B(z,{\overline \zeta})$ by the orthonormal series
$\sum \phi_\alpha (z) {\overline {\phi_\alpha(\zeta)}} $.  The restriction
of $K$ to $V_j$ is then given by
the kernel $f(z,{\overline \zeta}) ||\Phi^d(z)||^2 $. By the positivity and
by Proposition 2, 
it then follows that  $f(z,{\overline z}) ||\Phi^d(z)||^2$ is a
squared norm as in (21).  $ \spadesuit$

\bigskip

{\bf III. Interpretation in terms of line bundles}.

Let ${\bf P_{n}}$ denote the complex projective space of lines through the 
origin in ${\bf C^{n+1}}$. As usual, see [W] for example, for $j=0,1,...,n$ 
we let $U_j$ be the open subset where the homogeneous coordinate
$z_j$ is not zero. 
A line bundle $E$ over ${\bf P_n}$
is defined by its transition functions $g_{jk} : U_j \cap U_k \to {\bf C}^*$.
Recall that the universal line bundle
${\bf U}$ over ${\bf P_n}$ is the line bundle
whose transition functions are $ h_{jk}(z) ={z_j \over z_k}$.
 The hyperplane section bundle $H$ over 
${\bf P_n}$ is the dual line bundle; its 
transition functions are $ g_{jk} (z) = {z_k \over z_j}$. The $m$-th powers
of these bundles are the bundles ${\bf U}^m$ and $H^m$
whose transition functions are $({z_j \over z_k})^m$  and $({z_k \over z_j})^m$
respectively. The holomorphic sections of $H^m$ are homogeneous polynomials of degree $m$
in the homogeneous coordinates. 

A metric on a line bundle $E$ over ${\bf P_n}$ whose transition functions are $g_{jk}$ 
determines a positive 
function $p_j$ in each $U_j$ such that
$p_j = p_k |g_{kj}|^2 $ on the overlap $ U_j \cap U_k$, and 
conversely such a collection of positive
functions defines a metric on $E$. This is easily seen by observing that if $s_k$
denotes the representation of a local section in $U_k$, then its squared length
$p_k |s_k|^2$ must equal  $p_j |s_j|^2$ in the overlap $ U_j \cap U_k$. Since
$s_k = g_{kj} s_j$ we see that the metric must transform by

$$ p_j = p_k |g_{kj}|^2 $$

Let $ {\cal P}$ denote the set of bihomogeneous polynomials on ${\bf C^{n+1}}$ that
are positive away from the origin. 
When $f \in {\cal P}$,
we can use it to
define metrics on ${\bf U}^m$ or on $H^m$ over  ${\bf P_n}$ by the following method. 
For $j=0,...,n$ we put  
$$ f_j(z, {\overline z}) = { f(z,{\overline z}) \over |z_j|^{2m} } \eqno (U) $$
Of course $f_j$ is defined and positive in 
$U_j$. On $U_j \cap U_k$, we have 
$$ {f_j \over f_k} = |({z_k \over z_j})^m|^2 .$$
Therefore the collection of functions $f_j$ defines a metric on the line bundle whose
transition functions are $h_{jk} = ({z_k \over z_j})^m$. This is the bundle ${\bf
U}^m$.

Had we put 
$$ p_j(z, {\overline z}) = { |z_j|^{2m} \over f(z,{\overline z}) } \eqno (H) $$
then the functions $p_j$ would define a metric on the $m$-th power of the
hyperplane section bundle $H^m$.

{\bf Definition 1}. We say that a metric defined on ${\bf U}^m$ or $H^m$ is a 
{\it special metric} if it is defined by a bihomogeneous positive polynomial $f$ as
in (U) or (H) above. We write the bundle and the special metric on ${\bf U}^m$
as the ordered pair $({\bf U}^m, f)$.

Next we consider the effect of multiplying $f$ by the 
squared norm $||\Phi^d||^2$ in Theorem 2. Since $||\Phi^d||^2 \in {\cal P}$, 
it defines a metric on
${\bf U}^d$ or $H^d$ as above. 

In the case of the unit ball, we use $||z||^2$ to define a special metric on
${\bf U}$; this is the standard metric. 
When we raise the squared norm to the power $d$, we obtain metrics on
${\bf U}^d$ or $H^d$. This amounts to taking the $d$-fold tensor product of the bundle.
The resulting metric differs only by constants from the metrics defined by
$||\Phi^d||^2$.

We may consider $||z||^{2d} f(z,{\overline z}) $ as a metric on
${\bf U}^{m+d}$ or  $H^{m+d}$ in
the same way. Using tensor products we can write

$$ ({\bf U}^d, ||z||^{2d}) \otimes ({\bf U}^m, f(z,{\overline z})) = 
({\bf U}^{m+d}, ||z||^{2d} f(z,{\overline z})) $$

Suppose that $ ||z||^{2d} f(z,{\overline z}) =||g||^2 $ and the 
components of $g$ form a basis for $V_{m+d}$. 
Let $N=N(m+d,n+1)$
denote the dimension of
the vector space of homogeneous polynomials of 
degree ${m+d}$ in $n+1$ complex variables, and consider the universal bundle
${\bf U}$ over
${\bf P_N}$. The line bundle ${\bf U}^{m+d}$ over ${\bf P_n}$ is obviously the pullback
by $g$ of
 the line bundle ${\bf U}$ over 
${\bf P_N} $.  The metric is also given by a pullback.
If we equip ${\bf U}$ with the metric given by $||L(\zeta)||^2$, where $L$
is the appropriate invertible linear mapping, then
$({\bf U}^{m+d}, ||g(z)||^2)$ over 
${\bf P_n}$ is the pullback of $({\bf U},||L(\zeta)||^2)$
over ${\bf P_N}$.

We can now restate Theorem [CD] and also Theorem 2.

{\bf Theorem 3}. Let $({\bf U}^m, f)$ denote 
the $m$-th power of the universal line bundle over
${\bf P_n}$ with special metric defined by $f$.
Then there is an integer $d$ so that $({\bf U}^{m+d}, ||z||^{2d} f(z, {\overline z}))$
is a (holomorphic) pullback $g^*({\bf U}, ||L(\zeta)||^2)$
of the standard metric on the universal bundle over  ${\bf P_N}$. 
The mapping $g: {\bf P_n} \to {\bf P_N}$ is a holomorphic (polynomial) embedding and
$L$ is an invertible linear mapping.

$$ ({\bf U}^m, f) \otimes ({\bf U}^d, ||z||^{2d}) = 
({\bf U}^{m+d}, ||z||^{2d} f(z,{\overline z})) = ({\bf U}^{m+d}, ||g(z)||^2) $$
We have the bundles and metrics $$\pi_1 : ({\bf U}^m,f) \to {\bf P_n}$$
$$\pi_2 : ({\bf U}^{m+d},||z||^{2d} f) \to {\bf P_n}$$
$$\pi_3 : ({\bf U},||L(\zeta)||^2) \to {\bf P_N}$$

The first is not necessarily a pullback of the third, but for sufficiently large $d$,
the second must be. We conclude the paper by restating Theorem 2 in
this language.

{\bf Theorem 4}. Suppose that 
$\Omega$ is a smoothly bounded pseudoconvex
circled domain of finite type in ${\bf C^n}$.
Let  $({\bf U}^d, ||\Phi^d||^2)$ denote the $d$-th power 
of the universal bundle over ${\bf P_{n-1}}$ 
with special metric $||\Phi^d||^2$. Here 

$$ ||\Phi^d||^2 = \sum \phi_\alpha(z) {\overline {\phi_\alpha(\zeta)}} $$
and the sum is taken over an orthonormal basis for $V_d$.
Let $f$ define a special metric
on $({\bf U}^m,f) $. Then there is an integer $d_0$ such that, for all $d \ge d_0$,
the bundle $({\bf U}^{m+d}, f ||\Phi^d||^2)$ is the holomorphic
pullback of the universal line bundle $({\bf U},||L(\zeta)||^2)$ 
over some 
${\bf P_N}$ with the standard
metric. Again $L$ is an invertible linear mapping.

\bigskip

{\bf References}.

[C1] David Catlin, Subelliptic estimates for the ${\overline \partial}$-Neumann problem
on pseudoconvex domains, Annals of Math 126(1987), 131-191.

[C2] David Catlin, Necessary conditions for subellipticity of the ${\overline
\partial}$-Neumann problem, Annals of Math 117 (1983), 147-171.

[CD] David W. Catlin and John P. D'Angelo,
A stabilization theorem for Hermitian forms and applications to holomorphic
mappings, Math Research Letters 3 (1996), 149-166. 

[D] John P. D'Angelo, Several complex variables and the geometry of real hypersurfaces,
CRC Press, Boca Raton, 1993.

[De] Jean-Pierre Demailly, $L^2$ vanishing theorems for positive line bundles and
adjunction theory, CIME Session, Transcendental Methods in Algebraic Geometry, Cetraro
Italy, 1994.

[FK] Gerald B. Folland and J. J. Kohn, The Neumann problem for the
Cauchy-Riemann complex, Annals of Math Studies 75, Princeton University 
Press, 1972.

[HI] Gennadi M. Henkin and Andrei Iordan, Compactness of the Neumann operator for
hyperconvex domains with non-smooth boundaries, Math. Ann. 307 (1997), 151-169.

[Ke] Norberto Kerzman, The Bergman kernel: Differentiability at the boundary, Math
Annalen 195 (1972), 149-158.

[K] Joseph J. Kohn, Subellipticity of the ${\overline \partial}$-Neumann problem on
pseudoconvex domains: Sufficient conditions, Acta Mathematica, 142(1979), 79-122. 

[Ru] Walter Rudin, Functional Analysis, McGraw-Hill, New York, 1973.

[W] Raymond O. Wells, Differential Analysis on Complex Manifolds, Prentice-Hall,
Englewood Cliffs, New Jersey, 1973.

{\bf Author addresses}

Dept. of Mathematics, Purdue Univ., W. Lafayette IN 47907

catlin@math.purdue.edu

(Catlin supported by NSF grant DMS 94015480)

Dept. of Mathematics, Univ. of Illinois, Urbana IL 61801

jpda@math.uiuc.edu

(D'Angelo supported by NSF grant DMS 9304580 at IAS and by MSRI)

\end